\documentclass[graybox]{svmult}
\usepackage{type1cm}        
\usepackage{makeidx}         
\usepackage{graphicx}        
\usepackage{multicol}        
\usepackage[bottom]{footmisc}

\usepackage{newtxtext}       %
\usepackage{newtxmath}       

\makeindex             

\usepackage{amsbsy}
\usepackage{latexsym}
\usepackage[mathscr]{eucal}
\usepackage{amsfonts}
\usepackage{amssymb}
\usepackage{mathtools}

\graphicspath{{img/}}

\newcommand{\dw}[1]{\textcolor{magenta}{#1}}

\def\rnew{\color{black}}

\begin{document}
\bibliographystyle{plain}
\title*{
State Estimation - The Role of Reduced Models}
\author{Albert Cohen, Wolfgang Dahmen and Ron DeVore
\thanks{%
A.C. was supported by {\wnew ERC Adv Grant BREAD}. W.D.  was supported in part by the   NSF-grant DMS 17-20297,
by the Smart State Program of the State of South Carolina, and  the Williams-Hedberg Foundation. R.D. was supported by the NSF grant DMS 18-17603.  A portion of this work was completed when the authors  were supported as visitors to the Isaac Newton Institute of Cambridge University.
}  }

\institute{Albert Cohen \at Laboratoire Jacques-Louis Lions,, Sorbonne Universit\'{e}, 4, Place Jussieu, 75005 Paris, France, \email{cohen@ann.jussieu.fr}
\and Wolfgang Dahmen \at Mathematics Department, University of South Carolina, 1523 Greene Street, Columbia SC 29208, USA,\email{dahmen@math.sc.edu}
\and Ronald DeVore \at Department of Mathematics, Texas A \& M University, College Station, TX 77843-3368, USA,
\email{ronald.a.devore@gmail.com}
}

\hbadness=10000
\vbadness=10000
\newtheorem{prop}[lemma]{Proposition}
\newtheorem{cor}[lemma]{Corollary}
\newtheorem{proper}[lemma]{Properties}
\newtheorem{assumption}[lemma]{Assumption}
%
\newenvironment{disarray}{\everymath{\displaystyle\everymath{}}\array}{\endarray}

\def\vp{\varphi}
\def\<{\langle}
\def\>{\rangle}
\def\t{\tilde}
\def\i{\infty}
\def\e{\varepsilon}
\def\sm{\setminus}
\def\nl{\newline}
\def\o{\overline}
\def\wt{\widetilde}
\def\wh{\widehat}
\def\cT{{\cal T}}
\def\cA{{\cal A}}
\def\cI{{\cal I}}
\def\cV{{\cal V}}
\def\cB{{\cal B}}
\def\cF{{\cal F}}

\def\cR{{\cal R}}
\def\cD{{\cal D}}
\def\cP{{\cal P}}
\def\cJ{{\cal J}}
\def\cM{{\cal M}}
\def\cO{{\cal O}}
\def\Chi{\raise .3ex
\hbox{\large $\chi$}} \def\vp{\varphi}
\def\lsima{\hbox{\kern -.6em\raisebox{-1ex}{$~\stackrel{\textstyle<}{\sim}~$}}\kern -.4em}
\def\lsim{\hbox{\kern -.2em\raisebox{-1ex}{$~\stackrel{\textstyle<}{\sim}~$}}\kern -.2em}
\def\[{\Bigl [}
\def\]{\Bigr ]}
\def\({\Bigl (}
\def\){\Bigr )}
\def\[{\Bigl [}
\def\]{\Bigr ]}
\def\({\Bigl (}
\def\){\Bigr )}
\def\L{\pounds}
\def\pr{{\rm Prob}}
\newcommand{\cs}[1]{{\color{magenta}{#1}}}
\def\ds{\displaystyle}
\def\ev#1{\vec{#1}}     
\newcommand{\lt}{\ell^{2}(\nabla)}
\def\Supp#1{{\rm supp\,}{#1}}
\def\R{\mathbb{R}}
\def\E{\mathbb{E}}
\def\nl{\newline}
\def\T{{\relax\ifmmode I\!\!\hspace{-1pt}T\else$I\!\!\hspace{-1pt}T$\fi}}
\def\N{\mathbb{N}}
\def\Z{\mathbb{Z}}
\def\N{\mathbb{N}}
\def\Zd{\Z^d}
\def\Q{\mathbb{Q}}
\def\C{\mathbb{C}}
\def\Rd{\R^d}
\def\gsim{\mathrel{\raisebox{-4pt}{$\stackrel{\textstyle>}{\sim}$}}}
\def\sime{\raisebox{0ex}{$~\stackrel{\textstyle\sim}{=}~$}}
\def\lsim{\raisebox{-1ex}{$~\stackrel{\textstyle<}{\sim}~$}}
\def\div{\mbox{ div }}
\def\M{M}  \def\NN{N}                  
\def\L{{\ell}}               
\def\Le{{\ell^1}}            
\def\Lz{{\ell^2}}
\def\Let{{\tilde\ell^1}}     
\def\Lzt{{\tilde\ell^2}}
\def\Ltw{\ell^\tau^w(\nabla)}
\def\t#1{\tilde{#1}}
\def\la{\lambda}
\def\La{\Lambda}
\def\ga{\gamma}
\def\BV{{\rm BV}}
\def\Ga{\eta}
\def\al{\alpha}
\def\cZ{{\cal Z}}
\def\cA{{\cal A}}
\def\cU{{\cal U}}
\def\ms{{\rm ms}}
\def\wc{{\rm wc}}
\def\argmin{\mathop{\rm argmin}}
\def\argmax{\mathop{\rm argmax}}
\def\prob{\mathop{\rm prob}}
\def\A{\mathop{\rm Alg}}
\def\cen{{\rm cen}}

\def \bphi{{\bf\phi}}

\def\cO{{\cal O}}
\def\cA{{\cal A}}
\def\cC{{\cal C}}
\def\cS{{\cal F}}
\def\bu{{\bf u}}
\def\bz{{\bf z}}
\def\bZ{{\bf Z}}
\def\bI{{\bf I}}
\def\cE{{\cal E}}
\def\cD{{\cal D}}
\def\cG{{\cal G}}
\def\cI{{\cal I}}
\def\cJ{{\cal J}}
\def\cM{{\cal M}}
\def\cN{{\cal N}}
\def\cT{{\cal T}}
\def\cU{{\cal U}}
\def\cV{{\cal V}}
\def\cW{{\cal W}}
\def\cL{{\cal L}}
\def\cB{{\cal B}}
\def\cG{{\cal G}}
\def\cK{{\cal K}}
\def\cS{{\cal S}}
\def\cP{{\cal P}}
\def\cQ{{\cal Q}}
\def\cR{{\cal R}}
\def\cU{{\cal U}}
\def\bL{{\bf L}}
\def\bl{{\bf l}}
\def\bK{{\bf K}}
\def\bQ{{\bf Q}}
\def\bC{{\bf C}}
\def\X{X\in\{L,R\}}
\def\ph{{\varphi}}
\def\D{{\Delta}}
\def\H{{\cal H}}
\def\bM{{\bf M}}
\def\bx{{\bf x}}
\def\bj{{\bf j}}
\def\bG{{\bf G}}
\def\bS{{\bf S}}
\def\bP{{\bf P}}
\def\bW{{\bf W}}
\def\bX{{\bf X}}
\def\bT{{\bf T}}
\def\bV{{\bf V}}
\def\bv{{\bf v}}
\def\bt{{\bf t}}
\def\bz{{\bf z}}
\def\bw{{\bf w}}
\def \mvn {{\rm mvn}}
\def\rhom{{\rho^m}}
\def\diff{\hbox{\tiny $\Delta$}}
\def\EE{{\rm Exp}}
\def\lll{\langle}
\def\argmin{\mathop{\rm argmin}}
\def\argmax{\mathop{\rm argmax}}
\def\dJ{\nabla}
\newcommand{\ba}{{\bf a}}
\newcommand{\bb}{{\bf b}}
\newcommand{\bc}{{\bf c}}
\newcommand{\bd}{{\bf d}}
\newcommand{\bs}{{\bf s}}
\newcommand{\bff}{{\bf f}}
\newcommand{\bp}{{\bf p}}
\newcommand{\bg}{{\bf g}}
\newcommand{\by}{{\bf y}}
\newcommand{\br}{{\bf r}}
\newcommand{\be}{\begin{equation}}
\newcommand{\ee}{\end{equation}}
\newcommand{\bea}{$$ \begin{disarray}{lll}}
\newcommand{\eea}{\end{disarray} $$}
\def \Vol{\mathop{\rm  Vol}}
\def \mes{\mathop{\rm mes}}
\def\rad{\mathop{\rm rad}}
\def \Prob{\mathop{\rm  Prob}}
\def \exp{\mathop{\rm    exp}}
\def \sign{\mathop{\rm   sign}}
\newcommand{\mult}{\mathop{\rm   mult}}
\newcommand{\one}{\mathop{\rm   one}}
\def \wca{{\rm    wca}}
\def \msa{{\rm    msa}}
\def \sp{\mathop{\rm   span}}
\def \vphi{{\varphi}}
\def \csp{\overline \mathop{\rm   span}}
%
%
\newcommand{\beqn}{\begin{equation}}
\newcommand{\eeqn}{\end{equation}}
\def\beginproof{\noindent{\bf Proof:}~ }
\def\endproof{\hfill\rule{1.5mm}{1.5mm}\\[2mm]}

\newcommand{\utr}{u^{\rm true}}
\newcommand{\Cor}{\kappa}

\newcommand{\prox}{\mathrm{prox}}
\newcommand{\epi}{\text{epi}}
\newcommand{\pa}[1]{\left(#1\right)}
\newcommand{\xob}{x_{\mathrm{ob}}}
\newcommand{\xsol}{x^{\star}}
\newcommand{\xbar}{\bar{x}}
\newcommand{\xkm}{x_{k-1}}
\newcommand{\xk}{x_{k}}
\newcommand{\xkp}{x_{k+1}}
\newcommand{\xbarkp}{\bar{x}_{k+1}}
\newcommand{\tk}{t_{k}}
\newcommand{\tkp}{t_{k+1}}
\newcommand{\tbarkp}{\bar{t}_{k+1}}
\newcommand{\xtk}{\pa{\xk,\tk}}
\newcommand{\xtkp}{\pa{\xkp,\tkp}}
\newcommand{\xtbarkp}{\pa{\xbarkp,\tbarkp}}

\newcommand{\vsol}{v^{\star}}
\newcommand{\vbar}{\bar{v}}
\newcommand{\vkm}{v_{k-1}}
\newcommand{\vik}{v_{i,k}}
\newcommand{\vikp}{v_{i,k+1}}
\newcommand{\xiik}{\xi_{i,k}}
\newcommand{\xiikp}{\xi_{i,k+1}}
\newcommand{\vxik}{\pa{\vik,\xiik }_{i \in [p]}}
\newcommand{\vxiok}{\pa{\vik,\xiik}}
\newcommand{\vxikp}{\pa{\vikp,\xiikp}_{i \in [p]}}
\newcommand{\vxiokp}{\pa{\vikp,\xiikp}}

\newenvironment{Proof}{\noindent{\bf Proof:}\quad}{\endproof}

\renewcommand{\theequation}{\thesection.\arabic{equation}}
\renewcommand{\thefigure}{\thesection.\arabic{figure}}

\makeatletter
\@addtoreset{equation}{section}
\makeatother

\newcommand\abs[1]{\left|#1\right|}
\newcommand\clos{\mathop{\rm clos}\nolimits}
\newcommand\trunc{\mathop{\rm trunc}\nolimits}
\renewcommand\d{d}
\newcommand\dd{\mathrm d}
\newcommand\diag{\mathop{\rm diag}}
\newcommand\dist{\mathop{\rm dist}}
\newcommand\diam{\mathop{\rm diam}}
\newcommand\cond{\mathop{\rm cond}\nolimits}
\newcommand\eref[1]{{\rm (\ref{#1})}}
\newcommand{\iref}[1]{{\rm (\ref{#1})}}
\newcommand\Hnorm[1]{\norm{#1}_{H^s([0,1])}}
\def\int{\intop\limits}
\renewcommand\labelenumi{(\roman{enumi})}
\newcommand\lnorm[1]{\norm{#1}_{\ell^2(\Z)}}
\newcommand\Lnorm[1]{\norm{#1}_{L_2([0,1])}}
\newcommand\LR{{L_2(\R)}}
\newcommand\LRnorm[1]{\norm{#1}_\LR}
\newcommand\Matrix[2]{\hphantom{#1}_#2#1}
\newcommand\norm[1]{\left\|#1\right\|}
\newcommand\ogauss[1]{\left\lceil#1\right\rceil}
\newcommand{\QED}{\hfill
\raisebox{-2pt}{\rule{5.6pt}{8pt}\rule{4pt}{0pt}}%
  \smallskip\par}
\newcommand\Rscalar[1]{\scalar{#1}_\R}
\newcommand\scalar[1]{\left(#1\right)}
\newcommand\Scalar[1]{\scalar{#1}_{[0,1]}}
\newcommand\Span{\mathop{\rm span}}
\newcommand\supp{\mathop{\rm supp}}
\newcommand\ugauss[1]{\left\lfloor#1\right\rfloor}
\newcommand\with{\, : \,}
\newcommand\Null{{\bf 0}}
\newcommand\bA{{\bf A}}
\newcommand\bB{{\bf B}}
\newcommand\bR{{\bf R}}
\newcommand\bD{{\bf D}}
\newcommand\bE{{\bf E}}
\newcommand\bF{{\bf F}}
\newcommand\bH{{\bf H}}
\newcommand\bU{{\bf U}}
\newcommand\cH{{\cal H}}
\newcommand\sinc{{\rm sinc}}
\def\enorm#1{| \! | \! | #1 | \! | \! |}

\newcommand{\dm}{\frac{d-1}{d}}

\let\bm\bf
\newcommand{\balpha}{{\mbox{\boldmath$\alpha$}}}
\newcommand{\bbeta}{{\mbox{\boldmath$\beta$}}}
\newcommand{\bal}{{\mbox{\boldmath$\alpha$}}}
\newcommand{\bbi}{{\bm i}}

\newcommand{\cY}{\mathcal{Y}}
\newcommand{\U}{\mathbb{U}}
\newcommand{\V}{\mathbb{V}}
\newcommand{\W}{\mathbb{W}}

\newcommand{\test}{\text{test}}
\newcommand{\greedy}{\text{greedy}}

\def\nnew{\color{black}}
\def\mnew{\color{orange}}
\def\wnew{\color{black}}
\def\dw{\color{purple}}

\newcommand{\dI}{\Delta}
\maketitle
 
\abstract{
The exploration of complex physical or technological processes usually requires exploiting available information from different sources:
(i) {\em physical laws} often represented as a family of {\em parameter dependent partial differential equations} 
 and (ii)  {\em data} provided by {\em measurement devices} 
or {\em sensors}. The amount of sensors is typically limited and data acquisition may be expensive and in some cases even harmful.
This article reviews some recent developments for this ``small-data'' scenario where
inversion is strongly aggravated by the typically large parametric dimensionality. The proposed concepts may be viewed as 
exploring alternatives to {\em Bayesian   inversion} in favor of more deterministic accuracy quantification related to
the required computational complexity. We discuss {\it optimality criteria} which delineate intrinsic information limits,
and highlight the role of {\em reduced models} for developing efficient computational strategies. In particular, the need to {\em adapt} the reduced models - not to
a specific (possibly noisy) data set but rather to the sensor system -  is a central theme. This, in turn, is  facilitated by exploiting
{\em geometric} perspectives based on proper {\em stable variational formulations} of the continuous model.  
}

\section{Introduction}\label{sec:intro}
Modern sensor technology and data acquisition capabilities generate an ever  increasing
wealth of   data about  virtually every branch of science and social life.  Machine learning offers novel techniques for extracting quantifiable information from
such  large data sets.   While machine learning has already had  a transformative impact on a diversity of application areas in the ``big-data'' regime, particularly in image classification and artificial intelligence, it  is  yet to have  a similar impact  in many other areas of science.

 Utilizing data observations in the analysis of  scientific processes differs from  traditional learning in  that one has the additional information that these processes    are described by mathematical models - systems of partial differential equations (PDE)
or integral equations - that encode the  physical laws that  
govern     the  process. Such models, however, are often deficient, inaccurate, incomplete or  need to
be further calibrated by  determining a   large number of {\em parameters} in order to accurately represent an observed process.  
Typical guiding examples are Darcy's equation for the pressure in ground-water flow or electron impedance tomography.
Both are  based on second order elliptic equations as core models. The diffusion coefficients in these examples describe premeability or conductivity, respectively. The parametric representations of the coefficients
  could  arise, for instance, from Karhunen-Lo\`{e}ve expansions
of a random field that   represent ``unresolvable'' features to be captured by the model. In this case the number of parameters 
could actually be {\em infinite}. 

The use of  machine learning to describe  complex states of interest or even the
underlying laws, solely through data, seems to bear little hope. In fact, data acquisition is often expensive or even harmful 
as in applications involving radiation. Thus, a severe undersampling poses principal obstructions to {\em state} or {\em parameter estimation} by solely processing observational data through  
standard machine learning techniques.
It is therefore more  natural to try to effectively combine  the data information   with the  knowledge  of the underlying
physical laws   represented by parameter dependent families of PDEs. 

Methods   that  fuse  together data-driven and model-based approaches fall roughly into two categories. One  prototype of
a {\em data assimilation} scenario arises in meteorology where data are used to stabilize otherwise chaotic dynamical systems,
typically with the aid of (stochastic)  filtering techniques. 
A second setting,  in line with the above examples,
uses an underlying {\em stable} continuous model   to  {\em regularize} otherwise ill-posed estimation tasks in a 
``small-data'' scenario.  
{\em Baysian inversion} is a prominent   way of regularizing such problems.  It  relaxes  the estimation task to asking only for {\em posterior probabilities} of states or parameters to explain given observations.

 The present article reviews some recent developments on  data driven state and parameter estimation
that  can be
viewed as seeking alternatives to Bayesian inversion by placing  a  stronger focus on deterministic   uncertainty quantification
 and its relation to   {\em computational complexity}. 
The emphasis is on foundational  
aspects    such as the optimality of algorithms (formulated in an appropriate sense)  when treating   estimation tasks for ``small-data'' problems in {\em high-dimensional parameter} regimes. Central issues concern
the role of {\em reduced modeling} and the exploitation of intrinsic problem metrics provided by the {\em variational formulation}
of the underlying continuous family of PDEs. This is used by the so called   {\em Parametrized Background Data-Weak} (PBDW) framework, introduced in   \cite{MPPY}
and further 
analyzed in \cite{BCDDPW2},  to  identify  a suitable trial (Hilbert) space  $\U$  that accomodates the states and eventually also the data.  
An important point is to distinguish between
the {\em data} and corresponding {\em sensors} - here linear functionals in the dual $\U'$ of $\U$ -  from  which the data are generated.
This will be seen to actually open a {\em geometric} perspective   that sheds light on intrinsic estimation limits.
Moreover, in the deterministic setting, a pivotal role is played by the so called {\em solution manifold}, which is the set of all states that can be attained when the parameters in the PDE traverse the whole parameter domain.

Even with full knowledge   of  a state in the solution manifold,  to infer from it a corresponding
parameter is a {\em nonlinear severly ill-posed} problem typically formulated as a {\em non-convex} optimization problem. 
On the other hand, state estimation   from data is a {\em linear}, and hence a more benign inversion task mainly suffering under the current premises from a severe undersampling. We will, however, indicate how to reduce, under certain circumstances, the latter to the former problem  so as to end up with a {\em convex} optimization
problem. This motivates focusing in what follows mainly on state estimation. 
 A central question then becomes how to best invoke 
knowledge on the solution manifold to regularize the estimation problem without introducing unnecessarily ambiguous bias.
Our principal viewpoint is to recast state estimation as an {\em optimal recovery} problem which then naturally leads
one to explore the role and potential of   {\em reduced modeling}. 

The layout of the paper is as follows. Section \ref{sec:2} describes  the conceptual framework for {\em state estimation}
as an {\em optimal recovery task}.  This formulation allows the identification of  lower bounds for the best achievable recovery accuracy.

Section \ref{sec:one} reviews recent developments concerning a certain {\em affine} recovery scheme and highlights 
the role of {\em reduced models} adapted to the recovery task. The overarching theme is to establish  certified recovery bounds. 
 When striving for optimality of such affine recovery maps,  high parameter dimensionality is identified
as a major challenge.  We outline a recent remedy  that avoids the {\em Curse of Dimensionality} by trading
deterministic accuracy guarantees against analogs that hold with quantifiable high probability.

Even optimal affine reduced models can, in general, not be expected to realize the benchmarks identified in Section \ref{sec:2}.
To put the results in Section \ref{sec:one} in proper perspective, we comment in Section \ref{sec:pwa} on ongoing work
that uses the results on affine reduced models and corresponding estimators as a central building block 
  for nonlinear estimators. We also indicate briefly some ramifications
on parameter estimation.

\section{Models and data}\label{sec:2}
\subsection{The model}\label{ssec2.1}
Technological design or simulating physical processes is often based on  {continuum} models given   by a family
\be
\label{model}
  \cR(u,y)= 0, \quad y\in \cY,
\ee
  of partial differential Equations (PDEs) that depend on parameters $y$ ranging over a parameter domain $ \cY  \subset \R^{d_y}$.   We will always assume {\em uniform well-posedness}
of \eqref{model}:  for each $y\in \cY$, there exists  a unique solution $u=u(y)$  in some 
trial Hilbert space $\U$ which satisfies $\cR(u(y),y)= 0$. 

Specifically, we consider
only linear problems of the form 
$
\cB_y u=f, 
$   
that is, 
\be
\cR(u,y)= f- \cB_y u.
\ee
Here $f$ belongs to the dual $\V'$ of a suitable {\em test space} $\V$ and $\cB_y$ is
a linear operator acting from $\U$ to $\V'$ that depends on $y\in\cY$.  Here, uniform well-posedness means then that $\cB_y$ is boundedly invertible
with bounds independent of $y$. 
By the  Babu$\check{s}$ka-Banach-Ne$\check{c}$as Theorem,
this is equivalent to saying that the bilinear form  
\be
 \label{SVF}
(u,v)\mapsto b_y(u,v):= (\cB_y u)(v) 
 \ee 
 satisfies the following {\em continuity} and {\em inf-sup conditions}
 \be
 \label{Cinf-sup}
   \sup_{u\in\U}\sup_{v\in \V} \frac{b_y(u,v)}{\|u\|_\U\|v\|_\V} \le C_b\quad {\rm and} \quad \inf_{u\in\U}\sup_{v\in \V} \frac{b_y(u,v)}{\|u\|_\U\|v\|_\V}\ge c_b>0,\quad y\in\cY,
 \ee
 together with the property that $b_y(u,v)=0, u\in \U$, implies $v=0$ (injectivity of $\cB_y^*$). 
  The relevance of this stability notion lies in  the entailed validity of the 
   {\em error-residual relation}
  \beqn
 \label{4}
 C_b^{-1}\|f- \cB_y v\|_{\V'} \le \|u(y)-v\|_\U\le c_b^{-1}\|f-\cB_y v\|_{\V'},\quad v\in \U,\,y\in \cY,
 \eeqn
 where  $\|g\|_{\V'}:= \sup\{g(v)\, : \, \|v\|_\V= 1\}$. Thus, errors in the trial norm are equivalent to residuals in the dual test norm which will be exploited in what follows.

For a wide range of problems such as space-time variational formulations, e.g. of parabolic or convection-diffusion problems, indefinite or singularly perturbed problems, the identification of a suitable pair $\U,\V$ that guarantees stability in the above sense is not entirely straightforward. In particular,    trial and test space may have to differ from each other, see e.g. \cite{BS,CDW,DHSW,SW}  for examples as well as some general principles.

 The simplest example, used for illustration purposes, is the {\em elliptic} family 
\be
\label{Poisson}
  \cR(u,y) = f + \div(a(y)\nabla u),
\ee
 set in  $\Omega \subset \R^{d_x}$ where $d_x\in \{1,2,3\}$, with boundary conditions $u|_{\partial\Omega}=0$. 
Uniform well-posedness follows then for $\U=\V= H^1_0(\Omega)$
 if we have for some fixed constants $0<r\le R<\infty$    the bounds
\be
\label{rR}r\le a(x,y)\le R,\quad (x,y)\in \Omega\times \cY,
\ee
readily implying \eref{Cinf-sup}.  

Aside from well-posedness, a second important structural property of the model \eref{model} is {\em affine parameter dependence}.
By this we mean that 
\be
\label{Raff}
 \cB_y u= \cB_0 u + \sum_{j=1}^{d_y} y_j \cB_j u,\quad y= (y_j)_{j=1,\dots,d_y}\in \cY,
\ee
where the operators $\cB_j:\U \to \V'$ are {\em independent} of $y$. In turn, the residual has a similar
affine dependence structure
\be
\label{Raff1}
\cR(u,y)= \cR_0 (u) + \sum_{j=1}^{d_y} y_j \cR_j u,\quad \cR_0(u):=f-\cB_0 u,\quad \cR_j=-\cB_j.
\ee
For the example \eref{Poisson} such a structure is encountered  for
{\em affine} parametric representations of the diffusion coefficients 
\be
\label{affine}
a(x,y)= a_0(x) +{ \sum_{j=1}^{d_y}}  y_j\theta_j(x), \quad (x,y)\in \Omega \times \cY,
\ee
i.e., the field $a$ is expanded in terms of some given spatial basis functions $\theta_j$. 
As indicated earlier,  the pressure equation in Darcy's law for porous media flow is an example for \eref{Poisson} where the diffusion coefficient $a(y)$ of the 
form \eqref{affine} may arise from
a stochastic model for permeability via a Karhunen-Lo\`{e}ve expansion.
In this case (upon proper normalization) $y\in [-1,1]^\N$ has, in principle, {\em infinitely} many entries, 
  that is $d_y=\infty$.
However, due to \eref{rR}, the $\theta_j$ should then have some decay as $j\to \infty$ which means that the parameters become less and less 
important when $j$ increases. 
Another example is electron impedance tomography involving the same type of elliptic operator where parametric expansions represent possible variations of conductivity often modeled as piecewise constants, i.e.,  the $\theta_j$ could be 
characteristic functions subordinate to a partition of $\Omega$.
In this case data are acquired through  sensors that act through trace functionals greatly adding to ill-posedness.

A central role in the subsequent discussion is played by the solution manifold
\be
\label{M}
\cM=u(\cY) := \{u(y)\,: \, y\in \cY\}
\ee
which is then the range of the {\em parameter-to-solution map} $u: y\mapsto u(y)$ comprised of all states that can be attained when 
$y$ traverses $\cY$. 
Without further mention,  $\cM$ will be assumed to be compact which actually follows 
under  standard assumptions met in all above mentioned examples.

Estimating states in $\cM$ or corresponding parameters from measurements requires the efficient approximation of elements in $\cM$.
 A common challenge encountered in all such models lies in the inherent {\em high-dimensionality} of the states 
$u=u(\cdot,y)$ as functions of $d_x$ spatial variables $x\in\Omega$ and $d_y\gg 1$ parametric variables $y\in \cY$.
In particular, when $d_y=\infty$ any calculation, of course, has to work with finitely many ``activated'' parameters whose number, however,
has to be coordinated with the spatial resolution of a numerical scheme to retain {\em model-consistency}. It is especially this issue
that hinders standard approaches based on  {\em first dicretizing} the parametric model because rigorously balancing spatial and parametric
uncertainties becomes then difficult. 

What renders such problem scenarios nevertheless numerically tractable is  a further property that will be implicitly assumed 
in what follows, namely  that 
the {\em Kolmogorov $n$-widths} of the solution manifold   
\be
\label{widths}
d_n(\cM)_\U := \inf_{{\rm dim}\, \U_n =n}\sup_{u\in \cM}\inf_{v\in\U_n}\|u-v\|_\U 
\ee 
exhibits at least some algebraic decay
\be
\label{s}
d_n(\cM)_\U  \lesssim n^{-s}
\ee 
for some $s>0$, see  \cite{CD} for a comprehensive account. 

For instance, this is known to be the case for elliptic models \eref{Poisson} with
 \eqref{rR},   as a consequence of the results of sparse polynomial approximation of the parameter
 to solution map $y\mapsto u(y)$ established e.g. in \cite{CDS}. More generally, \iref{s} can be established
under a general holomorphy property of the parameter to solution map, as a consequence
of a similar algebraic decay assumed on the $n$-widths of the parameter set, see \cite{CDnw}. 
For a fixed finite number $d_y<\infty$ of parameters, under certain structural assumptions on the 
 parameter representations (e.g. piecewise constants on checkerboard partitions) one can even establish  (sub-) exponential  decay rates, see \cite{BC16} for more details. Assuming $s$ in \eref{s} to have a ``substantial'' size for any range of $d_y$, is therefore justified. 

In summary, the results discussed below are valid and practically feasible for  well posed linear models \eref{Cinf-sup}
with affine parameter dependence \eref{Raff1} whose solution manifolds have rapidly decaying $n$-widths \eref{s}.

\subsection{The data}\label{ssec2.2}

Suppose we are given data $\bw = (w_1,\ldots,w_m)^\top\in \R^m$ representing observations {\rnew of} an unknown state $u\in \U$
obtained through $m$ linearly independent linear functionals $\ell_i\in \U'$, i.e.,
\be
\label{li}
w_i = \ell_i(u),\quad i=1,\ldots,m.
\ee
Since in real applications data acquisition may be costly or harmful we assume that $m$ is {\em fixed}. 
The central task to be discussed in what follows is to recover from this information an estimate for the observed unknown
state $u$, based on the prior assumption that $u$ belongs to $\cM$ or is close to $\cM$.
Moreover,
to bring out the essence of this estimation task we assume   for the moment that the data are noise-free.


Following \cite{MPPY,BCDDPW2},
we first recast the data in a ``compliant'' metric, by introducing the Riesz representers $\psi_i\in \U$,
defined by 
$$
( \psi_i,v)_\U = \ell_i(v),\quad v\in \U, \quad i=1,\ldots, m,
$$
The   $\psi_i$ now span the $m$-dimensional subspace $\W \subset \U$ which we refer to as {\em measurement space},
and the information carried by the $\ell_i(u)$ is equivalent to that of the orthogonal projection $P_\W u$ of $u$ to $\W$. The decomposition
\be
\label{deco}
u = P_\W u + P_{\W^\perp} u,\quad u\in \U,
\ee 
thus contains a first term that is ``seen'' by the sensors and a second (infinite-dimensional) term which cannot be detected.  The decomposition \eref{deco} may 
be seen as a sensor-induced ``coordinate system''    thereby opening   up a {\em geometric perspective} that will prove very useful in what
follows.
State estimation can then be viewed as learning from samples $w:=P_\W u$ the unknown ``labels'' $P_{\W^\perp} u\in \W^\perp$.

In this article,   we are interested in how well we can approximate $u$ from the information that $u\in\cM$  and $P_\W u=w$ with $w$ given to us.   Any such approximation 
is given by a mapping $A:w\to A(w)\in \U$.   
The overall performance of recovery on all of $\cM$ by the mapping $A$ is typically measured in  the worst case setting, that is, 
 \be
\label{wcA}
E_{\wc }(A,\cM,\W)= \sup_{u\in\cM} \|u-A(P_\W u)\|_\U.  
\ee
The optimal recovery error on $\cM$ is then defined as
 \be
\label{wcE}
E_{\wc }(\cM,\W)= \inf_AE_{\wc }(A,\cM,\W),
\ee
where the infimum is over all possible recovery maps. Let us observe that  the construction of  recovery maps  can be restricted to be  of the form
\be
\label{Aform}
A :w\to A(w),\quad A(w)= w+ B(w),\quad \text{with}\,\, B:\W \to \W^\perp.
\ee
Indeed, given any recovery mapping $A$, we can write $A(w)= P_\W A(w)+P_{\W^\perp} A(w)$ and the performance of the recovery can only be improved if we replace the first term  by $w$. In other words, $A(w)$ should belong to the affine space 
\be
\U_w:=w+\W^\perp,
\ee
that contains $u$. The mappings $B$ are commonly referred to as liftings into $\W^\perp$. 
 
\subsection{Optimality criteria and numerical recovery}  
  Finding a best recovery map $A$ attaining \eref{wcE}  is  known as  {\it optimal recovery} .  The best mapping has a  well-known  simple theoretical description,
  see e.g. \cite{MR}, that we now describe.
Note first that
a precise recovery of the unknown state $u$ from the given information   is generally impossible.     Indeed, the best we can say about $u$ is that it lies in the {\it manifold slice} 
\be
\label{slice}
\cM_w:= \{u\in \cM : P_\W u =w\} = \cM\cap \U_w, 
\ee
which is  comprised of all elements in $\cM$ sharing the same measurement $w\in \W$. 
  The Chebyshev ball $B(\cM_w)$ is the smallest ball in $\U$ that contains $\cM_w$. The best recovery algorithm is then given by  the mapping
\be
\label{CC}
 A^*(w): =\cen(\cM_w),
\ee
that assigns to each $w\in \cM$ the center $\cen(\cM_w)$ of $B(\cM_w)$, called the Chebyshev center of $\cM_w$.  
Then, the radius ${\rm rad}(\cM_w)$ of $B(\cM_w)$ is the best worst case error over the class $\cM_w$
The best worst case error over $\cM$, which is achieved by $A^*$, is thus given by
\be
\label{d0}
E_{\wc}(\cM,W) 
:= \max_{w\in \W} {\rm rad}(\cM_w).  
\ee

 While the above mapping $A^*$ gives a nice theoretical description of the  optimal recovery algorithm, it is typically not numerically implementable since the Chebychev center $\cen(\cM_w)$ is not easily found. Moreover, such an optimal algorithm is  
highly nonlinear and possibly discontinuous.  The purpose of this section is to formulate
a more modest goal for the performance of a recovery algorithm with the hope that this more modest goal can be met with a numerically constructable algorithm.  The remaining sections of the paper introduce  numerically implementable recovery mappings, analyze their performance,  
and  evaluate  the numerical cost in constructing these mappings.

The search for a numerically realizable algorithm must out of necessity  lessen the performance criteria.   A first possibility is to  weaken the performance criteria  to  {\em near best}   algorithms.  {This  means that we search for an algorithm $A$ such that
\be
E_{\wc }(A,\cM,\W) \leq C_0E_{\wc }(\cM,\W),
\ee
with a reasonable value of $C_0>1$. 
For example, any mapping $A$ which takes $w$ into an element in the Chebyshev ball of  $\cM_w$ is near best with constant 
$C_0=2$. }
However,  finding near best mappings $A$ also  seems to be numerically out of reach.    

In order to formulate  a more attainable performance criterion, we return to our earlier observations about uncertainty in both the model class $\cM$ and  
in the measurements $w$.  The former is a modeling error while the latter is an inherent measurement error.
   Both of these uncertainties can be  quantified by introducing for each $\e>0$,   the $\e$-neighborhood of the manifold  
   \be
   \cM^\e := \{v\in \U: \dist\,(v,\cM)_\U\le \e\}.
\ee   
The uncertainty in the model can be thought of as saying the sought after $u$ is in $\cM^\e$ rather than $u\in \cM$.  Also, we may formulate uncertainty (noise) in the measurements as saying  that they are not measurements of  a $u\in\cM$ but rather some $u\in\cM^\e$.   Here the value of $\e$ quantifies these uncertainties.

 Our new  goal is to  numerically  construct  a recovery map   $A$ that  is near-optimal on $\cM^\e$,   for some given $\e>0$.  Let us note that $\cM^\e$ is not  compact. 
An algorithm $A$ is  worst-case near  optimal for  $\cM^\e$ if and only if its performance is bounded by a constant multiple of  the diameter 
\be
\label{deps}
\delta_\e(\cM,\W) := \max\,\{\|u-v\|_\U: u,v\in\cM^\e,\, P_\W(u-v)=0\}.
\ee
Notice  that $\e=0$ gives the performance criterion for near optimal recovery over $\cM$.  
One can show that the function $\e\mapsto \delta_\e(\cM,\W)$ is   monotone non-decreasing in $\e$, continuous from the right, and
$\lim_{\e\to 0^+} \delta_\e(\cM,\W)= \delta_0(\cM,\W)$.  
 The speed at which $\delta_\e(\cM,\W)$ approaches $\delta_0(\cM,\W)$ reflects the ``condition'' of the
estimation problem depending on $\cM$ and $\W$. While the practical realization of worst-case near-optimality for $\cM^\e$  is already
a challenge, quantifying corresponding computational cost would require assumptions on the condition of the problem.

One central theme, guiding subsequent discussions, is therefore to find recovery maps $A_\e$ that 
realize an error bound of the form
\be
\label{wish}
E_{wc}(A_\e,\cM,\W) \le  C_0 \delta_\e(\cM,\W).
\ee
Any a priori information on   
measurement accuracy and model bias might be used to choose a viable tolerance $\e$.  

High parametric dimensionality poses particular challenges to estimation tasks   when 
the targeted error bounds are in the above worst case sense. These challenges can be 
somewhat mitigated when adopting a Bayesian point of view \cite{Stuart}. 
The prior information on $u$ is then described by a probability distribution $p$ on $\U$,
which is supported on $\cM$. Such a measure is typically induced by a probability distribution on $\cY$ that may or may not
be known. In the latter case, sampling $\cM$, i.e., computing snapshots $u(y^i)$, $i=1,\ldots,N$, for i.i.d. samples $y^i\in\cY$,
provides labeled data $(w_i,w_i^\perp)= (P_\W u(y^i),P_{\W^\perp} u(y^i))$ according to  the sensor-based decomposition \eref{deco}.
This puts us into the setting of {\em regression} in machine learning asking for an estimator that predicts for any new measurement $w\in\W$ 
its lifting $w^\perp=B(w)$. It is then natural to measure the performance of an algorithm
in an averaged sense. The best estimator $A$ that minimizes the mean-square risk
\be
E_{\ms}(A,p,\W)=\E(\|u-A(P_\W u)\|^2)=\int_\U \|u-A(P_\W u)\|^2 dp(u)
\label{meansquare}
\ee
is given by the 
conditional expectation
\be
\label{exp}
A(w) = \E(u| P_\W u=w).
\ee
 Since always  $E_{\ms}(A,p,\W)\le E_\wc(A,\cM,\W)$, the optimality benchmarks are somewhat weaker. In the rest of this paper, we adhere to 
the worst case error in the deterministic setting that only assumes membership of $u$ to $\cM$ or $\cM^\e$.

The following section is concerned with an  important  {\em building block} on a pathway towards
achieving \eref{wish} at quantifiable computational cost. This building block, referred to as {\em one-space method} is a linear (affine)
scheme which is, in principle, simple and
  easy to numerically implement. It depends on suitably chosen subspaces. We highlight the {\em regularizing property}  of these subspaces 
as well as ways to {\em optimize} them. This will reveal certain intrinsic obstructions caused by 
{\em parameter dimensionality}.   The one-space method by itself will generally not achieve \eqref{wish}
but,   as indicated earlier, can be used as a building block in a {\em nonlinear} recovery scheme that may indeed meet the goal \eref{wish}.

\section{The one-space method}\label{sec:one}
 
 \subsection{Subspace regularization}\label{sec:3}

 The one space method can 
be {viewed}  as a simple regularizer for state estimation. The resulting recovery map is induced
   by an $n$-dimensional subspace $\U_n$ of 
$\U$ for $n\le m$.    Assume that,  for each $n\geq 0$, we are given a subspace $\U_n\subset \U$
 of dimension $n$ whose distance from $\cM$ can be assessed
 \be
 \label{3.1}
 \dist (\cM,\U_n)_\U:=  \max_{u\in \cM} \dist (u,\U_n)_\U\le \e_n.
 \ee
Then   the cylinder
 \be
 \label{Kn}
 \cK(\U_n,\e_n) := \{ u\in \U: \dist (u,\U_n)_\U \le \e_n\} 
 \ee
contains $\cM$ and likewise the cylinder $\cK(\U_n,\e_n+\e)$ contains $\cM^\e$.
Our prior assumption that the observed state  belongs to $\cM$ or $\cM^\e$ can then be relaxed by
 assuming membership to these larger but simpler sets. 
 
 Remarkably, one can now  realize  an optimal recovery map quite easily  
 that  meets the relaxed benchmark $E_\wc(\cK(\U_n,\e_n),\W)$: in \cite{BCDDPW2}
 it was shown that the Chebychev center of the slice 
 \be
 \cK_w(\\U_n,\e_n):=\cK(\U_n,\e_n)\cap \U_w,
 \ee 
 is exactly given by
 the state in $\U_w$ that is closest to $\U_n$, that is
\be
\label{u*}
u^*= u^*(w) := \argmin_{u\in\U_w}\|u - P_{\U_n}u\|_\U.
\ee
This minimizer exists and can be shown to be unique as long as $\U_n\cap \W^\perp=\{0\}$. The  
corresponding optimal recovery map
\be
\label{rep}
A_{\U_n} : w \mapsto u^*(w)
\ee
was first introduced in \cite{MPPY} as the Parametrized Background Data Weak (PBDW) algorithm,
and is referred to as the {\em one-space} method in \cite{BCDDPW2}.
Due to its above minimizing property, it is readily checked that this map is linear
and can be determined with the aid of the singular value decomposition of the cross-Gramian
between any pair of orthonormal basis for $\U_n$ and $\W$. 

The worst case error $E_\wc(\cK(\U_n,\e_n),\W)$ can be described more precisely by introducing 
\be
\label{mu}
\mu (\U_n,\W) := \sup_{v\in\U_n}\frac{\|v\|_\U}{\|P_\W v\|_\U}
\ee
which is finite if and only if $\U_n\cap \W^\perp=\{0\}$. This quantity, also introduced in a related but slightly different context in \cite{Adcock},  
is therefore related to the angle between the spaces $\U_n$ and $\W$. It becomes 
large when $\U_n$ contains elements that are nearly perpendicular to $\W$.
It is  actually computable: one has $\mu(\U_n,\W)= \beta(\U_n,\W)^{-1}$ where 
\be
\label{infsup}
\beta(\U_n,\W) := \inf_{w\in\W}\sup_{v\in\U_n} \frac{\langle v,w\rangle_\U}{\|v\|_\U \|w\|_\U},
\ee
and $\beta(\U_n,\W)$ is the smallest singular value of the cross-Gramian   between any pair of orthonormal bases for $\W$ and 
$\U_n$. It has been shown in \cite{BCDDPW2,MPPY} that the worst case error bound over 
$\cK(\U_n,\e_n)$
is given by
\be
\label{erroreps}
E_\wc(A_{\U_n},\cK(\U_n,\e_n),\W) =E_\wc(\cK(\U_n,\e_n),\W)= \mu (\U_n,\W)\e_n.
\ee
The quantity $\mu(\U_n,\W)$ also coincides with the norm of the linear recovery map $A_{\U_n}$.
Relaxing the prior $u\in\cM$ by exploiting information on $\cM$ solely through approximability of $\cM$ by $\U_n$, 
thus implicitly {\em regularizes} the estimation task:
whenever $\mu(\U_n,\W)$ is finite, the optimal recovery map $A_{\U_n}$ is bounded 
and hence Lipschitz. 

One important observation is that the map $A_{\U_n}$ is actually independent of $\e_n$. In particular
it achieves optimality for the smallest possible containement cylinder 
\be
\cK(\U_n): = \cK(\U_n,\dist(\cM,\U_n)_\U),
\ee
and therefore, since $E_\wc(A_{\U_n},\cM,\W)\leq  E_\wc(A_{\U_n},\cK(\U_n),\W) =E_\wc(\cK(\U_n),\W)$,
\be
E_\wc(A_{\U_n},\cM,\W)\leq  
\mu (\U_n,\W)\dist\,(\cM,\U_n)_\U.
\label{errorM}
\ee
Likewise, the containment $\cM^\e \subset \cK(\U_n,\dist\,(\cM,\U_n)_\U+\e)$ implies that
\be
E_\wc(A_{\U_n},\cM^\e,\W) \leq \mu (\U_n,\W)(\dist\,(\cM,\U_n)_\U+\e).
\ee

On the other hand, the recovery map $A_{\U_n}$ may be far from optimal over the sets $\cM$ or $\cM^\e$. This 
is due to the fact that the cylinders $\cK(\U_n,\e_n)$ and $\cK(\U_n,\e_n+\e)$ may be much larger
than $\cM$ or $\cM^\e$. In particular, it is quite possible that for a particular observation $w$, one has
${\rm rad} (\cM_w)\ll {\rm rad} (\cK_w(\U_n,\e_n))$.
Therefore, we cannot generally expect that the one space method achieves our goal \eref{wish}.
In particular, the condition $n\le m$, which is necessary to avoid that $\mu (\U_n,\W)=\infty$, limits the 
dimension of an approximating subspace $\U_n$
and therefore $\e_n$ itself is inherently  bounded from below. The ``dimension budget'' has therefore to be used wisely in order 
to obtain good performance bounds. This typically rules   out ``generic approximation spaces'' such as finite element spaces,
and raises the question which subspace $\U_n$ yields the best estimator when applying the above method. 

 %
 \subsection{Optimal affine recovery}\label{sec:BA}

  The results of the previous section bring forward the question as to what is the best choice of the space 
 $\U_n$ for the given $\cM$. On the one hand, proximity to $\cM$ is desirable since $\dist\,(\cM,\U_n)_\U$ enters the error bound. However,
 favoring proximity, may increase $\mu(\U_n,\W)$. 
 Before addressing this question systematically, it is important to note that the above results carry over verbatim 
 when $\U_n$ is replaced by an {\em affine space} $\U_n = \bar u + \wt\U_n$ where $\wt \U_n\subset \U$ is a linear 
 space. This means the reduced model $\cK(\U_n,\e_n)$   is of the form 
 $$
\cK(\U_n,\e_n):= \bar u + \cK(\wt\U_n,\e_n).
 $$
 The best worst-case recovery bound  is now given by
  \be
 \label{still}
 E_\wc( \cK(\U_n,\e_n),\W)= \mu(\wt\U_n,\W)\e_n.
 \ee
 Intuitively, this may help to better control the angle between $\W$ and $\U_n$ by 
 anchoring the affine space at a suitable location (typically near or on $\cM$). More importantly, it    helps in {\em localizing}
 models via parameter domain decompositions that will be discussed  later.  
 
 The one-space algorithm discussed in the previous section confines the ``dimensionality'' budget of the approximation
 spaces $\U_n$ to $n\le m$.   In view of \eqref{errorM}, to obtain an overall good estimation accuracy,
 this space can clearly not be chosen arbitrarily but should be well adapted 
 both to the solution manifold $\cM$ and to measurement space $W$, that is, to the given observation functionals giving rise to the data.
 
  A simple way of {\em adapting} a recovery space to $\W$ is as 
follows: suppose for a moment that we were able to construct for  $n=1,\dots,m$, 
a hierarchy of spaces 
$
\U^{\rm nb}_1\subset \U^{\rm nb}_2\subset \cdots \subset\U^{\rm nb}_m,
$ 
that approximate $\cM$  in a {\em near-best} way, namely
 \be
 \label{KM}
 \dist\,(\cM,\U_n^{\rm nb})_\U \le C d_n(\cM)_\U.
 \ee
We may compute along the way the quantities $\mu(\U^{\rm nb}_j,\W)$,
 then choose 
\be
\label{nstar}
n^* = \argmin_{n\le m} \mu(\U^{\rm nb}_n,\W)\dist\,(\cM,\U^{\rm nb}_n)_\U,
\ee
and take the map $A_{\U^{\rm nb}_{n^*}}$. 
We sometimes refer to this choice as {\em ``poor man's algorithm''}. It is not clear though whether $\U^{\rm nb}_{n^*}$ is indeed a 
near-best choice for state recovery by the one-space method. In other words, one may question whether 
\be
\label{Kbest}
E_{\wc} (A_{\U^{\rm nb}_{n^*}},\cM,\W) \le C \inf_{{\rm dim}\wt\U\le m}E_{\wc}(A_{\wt\U},\cM,\W),
\ee
holds with a uniform constant $C<\infty$.   In fact,  numerical tests strongly suggest otherwise, which motivated in  \cite{CDDFMN}
the following alternative to the poor man's algorithm.

Recall that a given linear space $\U_n$ determines the
linear recovery map $A_{\U_n}$. Likewise a given affine space $\U_n$ determines an affine recovery map $A_{\U_n}$. 
Conversely, it can be checked that an affine recovery map $A$ determines an affine space $\U_n$ that allows 
one to interpret the recovery schemes as a one-space method in the sense that $A=A_{\U_n}$.
Denoting by  $\cA$ the class of all affine mappings of the 
form
\be
\label{affscheme}
A (w) =  w+ z +Bw,
\ee
where $z\in \W^\perp$ and $B\in \cL(\W,\W^\perp)$ is linear, we might thus
as well directly look for a mapping that
minimizes    
\be
E_{\wc}(A,\cM,\W):= \sup_{u\in\cM}\|u- A(P_\W u)\|_\U = \sup_{u\in\cM} \|P_{\W^\perp} u - { z}- { B} P_\W u\|_\U=:{\cE(z,B)}
\label{direct}
\ee
over $\cA$, i.e., over all  ${ (z,B)}\in  \W^\perp \times \cL(\W,\W^\perp)$.  It can be shown that indeed a minimizing pair $(z^*,B^*)$ exists,
 i.e.,
$$
\cE(z^*,B^*) = \min_{A\in \cA} E_\wc(A,\cM,\W) =: E_{\wc,\cA}(\cM,\W).
$$
However, the minimization of   $E_{\wc}(A,\cM,\W)$ over $(z,B)\in \W^\perp\times \cL(\W,\W^\perp)$ is far from 
practically feasible. In fact, each evaluation of   $E_{\wc}(A,\cM,\W)$ requires exploring $\cM$ and $B$ can have a range in the infinite dimensional
space $\W^\perp$. In order to arrive at a computationally tractable problem,  one needs to
\begin{enumerate}
\item[(i)]
Replace $\cM$ by a finite set $\wt \cM\subset \cM$, that should be sufficiently dense. Denseness can be quantified by requiring that 
$\wt\cM$ is a $\delta$-net for $\cM$ for some $\delta>0$, i.e.,  for any $u\in \cM$, there exists $\t u\in\wt \cM$ such that $\|u-\t u\|_\U\leq \delta$.
\item[(ii)]
Choose a finite dimensional space $\U_L\subset \U$ that approximates $\cM$ to a desired precision
$\dist\,(\cM,\U_L)_\U\leq \eta$, and replace $\W^\perp$ by the finite dimensional complement 
\be
\label{tW}
\wt\W^\perp := (\U_L +\W)\ominus \W
\ee
of $\W$ in $\U_L +\W$.
\end{enumerate}
The resulting optimization problem
\be
\label{pract}
(\tilde z,\wt B) = \argmin_{(z,B)\in \wt\W^\perp \times \cL(\W,\wt\W^\perp)} \sup_{u\in \wt\cM^\delta} \|P_{\W^\perp} u - { z}- { B} P_\W u\|_\U.
\ee
can be solved by primal-dual splitting methods providing a $O(1/k)$ convergence rate. 

Due to the perturbations (i) and (ii)  of the ideal minimization problem, the resulting $(\tilde z,\wt B)$ is no longer optimal.
However, one can show that
\be
\label{general}
E_\wc(\wt A,\cM,\W) \le E_{\wc,\cA}(\cM,\W) + \eta +C\delta,
\ee
where the constant $C$ is the operator norm of $B$ minimizing  \iref{direct}.
On the other hand,  since the range of any affine mapping $A$ is an affine space of
dimension at most $m$, therefore contained in a linear space of dimension at most $m+1$, one always has
$E_{\wc,\cA}(\cM,\W)\ge d_{m+1}(\cM)_\U$. Therefore $(\tilde z,\wt B)$ satisfies a near-optimal bound
\be
\label{na}
{  E_\wc(\wt A,\cM,\W)} \lsim E_{\wc,\cA}(\cM,\W),
\ee
whenever $\eta$ and $\delta$ are picked such that
\be
\label{if}
 \eta\lsim d_{m+1}(\cM)_\U,\quad {\rm and} \quad \delta \lsim d_{m+1}(\cM)_\U.
\ee
The numerical tests in \cite{CDDFMN} for a model problem of the type \eref{Poisson} with piecewise constant checkerboard
diffusion coefficients and $d_y$ up to $d_y=64$ show that this recovery map exhibits significantly better accuracy
than the method based on \eref{nstar}. It even yields  smaller error bounds than the affine mean square estimator \eref{meansquare}.
The following section   discusses the numerical cost entailed by conditions like \eref{if}.

\subsection{Rate-optimal reduced bases}\label{ssec:rbm}
To keep the dimension $L$ of the space $\U_L$ in \eref{tW}
small, a near-best subspace $\U_L^{{\rm nb}}$ in the sense of \eqref{KM} would be highly desirable. Likewise the poor man's scheme
\eref{nstar} would benefit from  such subspaces.  Unfortunately, such near-best subspaces are not practically 
accessible.  The {\em reduced basis method} aims to construct subspaces which come close to near-optimality
in a sense that we further explain next. The main idea is to generate theses subspaces by a sequence of
elements picked in the manifold $\cM$ itself, by means of a {\it weak-greedy algorithm}
introduced and studied in \cite{Maday}. In an idealized form, this algorithm proceeds as follows: 
given a current space $\U^{\rm wg}_n={\rm span}\{u_1,\dots,u_n\}$, one takes 
$u_{n+1}= u(y_{n+1})$ such that, for some fixed  $\gamma \in ]0,1]$,
$ 
\|u_{n+1}- P_{\U_n}u_{n+1}\|_\U\ge \gamma \max_{u\in \cM}\|u-P_{\U_n}u\|_\U,
$ 
or equivalently 
\be
\label{weakgreedyy}
\|u(y_{n+1})- P_{\U_n}u(y_{n+1})\|_\U\ge \gamma \max_{y\in \cY}\|u(y)-P_{\U_n}u(y)\|_\U,
\ee
Then, one defines $\U^{\rm wg}_{n+1}={\rm span}\{u_1,\dots,u_{n+1}\}$.
While unfortunately, the weak greedy algorithm  does  in general not produce spaces satisfying  \eref{KM}, it  does come close. Namely,  
it has been shown in \cite{BCDDPW1,DePW} that the spaces $\U^{\rm wg}_n$ are {\em rate-optimal} 
in the following sense:
\begin{enumerate}
\item
For any $s>0$ one has
\be
d_n(\cM)_\U \le C (n+1)^{-s}, \;  n\geq 0    \implies \dist\,(\cM,\U^{\rm wg}_n)_\U\le  \wt C (n+1)^{-s}, \; n\geq 0,
\label{rpol}
\ee
where $\wt C$ depends on $C,s,\beta,\gamma$. 
\item
For any $\beta>0$, one has
\be
d_n(\cM)_\U \le C e^{-c n^\beta}, \;  n\geq 0    \implies \dist\,(\cM,\U^{\rm wg}_n)_\U\le  \wt C e^{-\tilde c n^\beta}, \; n\geq 0,
\label{rexp}
\ee
where the constants $\t c, \wt C$ depend on $c,C,\beta,\gamma$.
\end{enumerate}
In the form  described above, the weak-greedy concept seems infeasible since it would, in principle, require  computing
the solution $u(y)$ for all values of $y\in\cY$ exactly,    exploring the whole exact solution manifold.
  However, its practical applicability is facilitated when there exists a  {\em tight} surrogate 
$R(y,\U_n)$, satisfying
\be
\label{surr}
c_R R(y,\U_n) \le \|u(y)-P_{\U_n}u(y)\|_\U=\dist\,(u(y),\U_n)\le c_R R(y,\U_n),\quad y\in \cY,
\ee
for uniform constants $0<c_R\le C_R<\infty$, which can be evaluated at affordable cost. Then, maximization of $R(y,\U_n)$ over $\cY$
amounts to the weak-greedy step \eref{weakgreedyy} with $\beta:=\frac {c_R}{C_R}$.
  According to \cite{DPW}, the validity of the following two conditions indeed allows one to derive 
computable  surrogates that satisfy \eref{surr}:
\begin{enumerate}
\item
The underlying parametric family of PDEs \eref{model} permits a uniformly {\em stable variational formulation} \eref{Cinf-sup}, 
and one has {\em affine parameter dependence} \eref{Raff1};
\item
The discrete projection  $\Pi_{\U_n}$ (of Galerkin or Petrov-Galerkin type) has the {\em best approximation property}, i.e.,  resulting errors
are uniformly comparable to the best approximation error. 
\end{enumerate}
Conditions (i) and (ii) ensure, in view of \eref{4}, that   $\|u(y)- P_{\U_n}u(y)\|_\U \sim \|\cR(y,\Pi_{\U_n}u(y))\|_{\V'}$ holds 
uniformly in $y\in\cY$.  Thus,
\be
\label{resup}
R(y,\U_n) := \|\cR(y,\Pi_{\U_n}u(y))\|_{\V'} = \sup_{v\in\V} \frac{\cR(y,\Pi_{\U_n}u(y)) (v)}{\|v\|_\V}
\ee
satisfies \eref{surr} and is therefore  a tight surrogate for $\dist\,(\cM,\U_n)_\U$. 
In the elliptic case \eref{Poisson}
under assumption \eref{rR}, then (i) and (ii) hold and the above comments reflect standard practice. 
For the wider scope of stable but {\em unsymmetric} variational formulations \cite{D,SW,BS}
the inf-sup conditions \eref{Cinf-sup} imply (i), 
but the Galerkin projection in (ii) needs to be replaced by a stable {\em Petrov-Galerkin} projection with respect to suitable test spaces $\V_n$
accompanying the reduced trial spaces $\U_n$. It has been shown in \cite{DPW} how to generate such test spaces 
with the aid of  a {\em double-greedy} strategy, see also \cite{D}.

  The main pay-off of using the surrogate $R(y,\U_n)$ is that one no longer needs to compute $u(y)$ but
only the low-dimensional projection $\Pi_{\U_n}u(y)$ by solving for each $y$ an $n\times n$ system,
which itself can be   rapidly assembled thanks to the affine parameter dependence \cite{Patera}.
However, one still faces the problem of
its exact maximization over $y\in\cY$. A standard approach is to  maximize instead
over a {\em discrete} training set $\wt\cY_{n}\subset Y$, which in turn induces a discretization of the
solution manifold
\be
\wt\cM_n=\{u(y)\, :\, y\in \wt\cY_n\}.
\ee
The resulting weak-greedy algorithm can be shown to remain rate optimal in the
sense of \iref{rpol} and \iref{rexp} if the discretization is fine enough so that 
$\wt\cM_n$ constitutes an $\e_n$-approximation net of $\cM$ where
$\e_n$ does not exceed $c\dist\,(\cM,\U^{\rm wg}_n)_\U$ for a suitable constant $0<c<1$.
In the current regime of {\em large or even infinite parameter dimensionality},
this becomes prohibitive because $\#\wt\cY_n$ would then typically scale like 
$O\big(\e_n^{-cd_y}\big)$, \cite{CDD-rand}.

As a remedy it has been proposed in   \cite{CDD-rand}
  to use training sets $\wt\cY_n$ that are generated
by {\em randomly sampling} $\cY$, and ask that the objective of rate optimality 
is met with high probability. This turns out to be achievable 
with training sets of much less prohibitive size.
In an informal and simplified manner the main result can be stated as follows.

\begin{theorem}
\label{th:random}
Given any target accuracy $\e >0$ and some $0<\eta<1$, then the weak greedy reduced basis algorithm based on choosing
at each step $N = N(\e,\eta)\sim |\ln \eta|+|\ln \e|$ randomly chosen training points in $\cY$ has the following properties with probability at least $1-\eta$: it terminates with $\dist\,(\cM,\U_{n(\e)})_\U \le {\rnew \e}$ as soon as the maximum of the surrogate over the current training set
falls below $c\e^{1+a}$ for some $c,a>0$. Moreover, if $d_n(\cM)_\U \le Cn^{-s}$, then $n(\e)\lsim \e^{-\frac 1s -b}$. The constants
$c,a,b$ depend on the constants in \eref{surr}, as well as on the rate $r$ of polynomial approximability of the
parameter to solution map $y\mapsto u(y)$. The larger $s$ and $r$, the smaller $a$ and $b$, and the closer the performance becomes to the ideal
one.
\end{theorem}

 \section{Nonlinear models}\label{sec:pwa}
 
 \subsection{Piecewise affine reduced models}

As already noted, schemes based on linear or affine reduced models of the form $\cK(\U_n,\e)$ can, in general, not
 be expected to realize the benchmark \eref{wish}, discussed earlier in Section \ref{sec:2}.    The 
convexity of the containment set $\cK(\U_n,\e)$ may cause the reconstruction error to be significantly larger than
$\delta_\e(\cM,\W)$.  Another way of understanding this limitation is that in order to make $\e$ small, one is enforced to
raise the dimension $n$ of $\U_n$, making the quantity $\mu(\U_n,\W)$ larger and eventually infinite
if $n>m$.

To overcome this principal limitation one needs to resort to {\em nonlinear} models
that better capture the non-convex geometry of $\cM$. 
One natural approach consists in replacing  
the single space $\U_n$ by a family $(\U^k)_{k=1,\dots,K}$ of affine spaces
\be
\U^k=\o u_k+\wt\U^k, \quad \dim(\wt\U^k)=n_k \leq m,
\ee
each of which aims to approximate a {\em portion} $\cM_k$ of $\cM$ to a prescribed target
accuracy simultaneously controlling  $\mu(\U^k,\W)$: fixing $\e>0$, 
we assume that we have at hand a partition of $\cM$ into portions 
\be
\label{Mi}
\cM = \bigcup_{k=1}^{K} \cM_k
\ee
such that 
\be
\label{Ui}
\dist\,(\cM_k,\U^k)_\U \le \e_k, \quad{\rm and}\quad  \mu(\wt \U^k,\W)\e_k \leq \e,  \quad k=1,\ldots,K.
\ee
One way of obtaining such a partition is trough a greedy splitting procedure
of the domain $\cY= [-1,1]^{d_y}$ which is detailed in \cite{CMN}.
The procedure terminates when for each cell $\cY_k$
the corresponding portion of the manifold $\cM_k$ can be associated 
to an affine $\U_k$ satisfying these properties.  We are ensured that this eventually occurs
since for a sufficiently fine cell $\cY_k$ one has ${\rm rad}(\cM_k)\leq \e$ which means
that we could then use a zero dimensional affine space $\U_k=\{\bar u_k\}$ 
for which we know that $\mu(\wt \U^k,\W)=1$.  
In this piecewise affine model, the containement property is now
\be
\cM\subset \bigcup_{k=1}^K \cK(\U_k,\e_k).
\ee
and the cardinality $K$ of the
partition depends on the prescribed $\e$. 

For a given measurement $w\in \W$, we may now
compute the state estimates
\be
u^*_k(w)= A_{\U^k}(w), \quad k=1,\dots,K,
\ee 
by the affine variant of the one-space method from \eref{u*}. Since $u\in \cM_{k_0}$ for some value $k_0$, we are ensured
that 
\be
\|u-u^*_{k_0}(w)\|_\U\leq \e,
\ee
for this particular choice. However $k_0$ is unknown to us and one has to rely on the 
data $w$ in order to decide which one among the affine models  is most appropriate for the
recovery. One natural {\em model selection} criterion can be derived if 
for any $\o u\in \U$ we have at our disposal a 
computable surrogate $S(\o u)$ that is equivalent to the distance from $\o u$
to $\cM$, that is
\be
\label{surr6}
c S(\bar u) \le \dist\,(\bar u,\cM)_\U\le C S(\bar u), \quad \dist\,(\bar u,\cM)_\U=\min_{y\in \cY} \|\o u-u(y)\|_\U,
\ee
for some fixed $0<c\leq C$. We give an instance of such a computable surrogate in \S \ref{ssec:mpr} below.
The selection criterion then consists in picking $k^*$ minimizing this surrogate between the
different available state estimates, that is, 
\be
\label{u*pwa}
u^*(w) := u^*_{k^*}(w)=\argmin\,\{S(u_k^*(w)) \,:\, k=1,\dots,K\}.
\ee
The following result,  established in \cite{CMN},  shows that this estimator 
now realizes the benchmark \eref{wish}
up to a multiplication of $\e$ by $\kappa := C/c$, where $c,C$ are the constants from \eref{surr6}.
%
\begin{theorem}
\label{th:wish}
Assume that \eref{Mi} and \eref{Ui} hold.
For any $u\in\cM$, if $w= P_\W u$, one has
\be
\label{wishtrue}
\|u- u^*(w)\|\le \delta_{\kappa\e}(\cM,\W),
\ee
where $\delta_\e(\cM,\W)$ is given by \eref{deps}.
 \end{theorem}

\subsection{Approximate metric projection and parameter estimation}\label{ssec:mpr}

A practically affordable realization of the surrogate $S(\o u)$ providing a  
{\em near-metric projection distance} to $\cM$ is a key ingredient of the
above nonlinear recovery scheme. Since it has further useful implications we add a few comments on that matter. 

As already observed in \eref{4}, whenever \eref{model} admits a stable variational formulation with respect to a suitable
pair $(\U,\V)$ of trial and test spaces,  the distance of any $\o u\in \U$
to any $u(y)\in \cM$ is uniformly equivalent 
to the residual of the PDE in $\V'$
\be
\label{err-res}
\|u(y)-\bar u\|_\U \sim \|\cR(\bar u,y)\|_{\V'},
\ee
Assume in addition that $\cR(u,y)$ depends {\em affinely} on $y\in\cY$, according to \eref{Raff1}.
Then, minimizing  $ \|\cR(\bar u,y)\|_{\V'}$ over $y$ is equivalent to solving a {\em constrained least squares}  problem 
\be
\label{LS}
\bar y = \argmin_{y\in\cY}\|\bg- \bM y\|_2,
\ee 
where $\bM$ is  a matrix of size $d_y\times d_y$
resulting from Riesz-lifts of the functionals $\cR_j (\bar u)$. 

The solution to this problem therefore satisfies \be
 \label{nearbest}
 \|\bar u - u(\bar y)\|_\U \le \kappa \inf_{y\in \cY} \|\bar u - u(y)\|_\U=\kappa \dist\,(\bar u,\cM)_\U. 
 \ee
where $\kappa$ is the quotient between the equivalence constants in \eref{err-res}.
The computable surrogate for the metric projection distance
of $\o u$ onto $\cM$ is therfore provided by
\be
S(\bar u)= \|\bar u - u(\bar y)\|_\U.
\ee

Since solving the above problem provides   an admissible parameter value $\o y\in \cY$,
this also has some immediate bearing on {\em parameter estimation}.
Suppose we wish to
estimate from $w= P_\W u(y^*)$ the unknown parameter $y^*\in \cY$. Assume further
that $A$ is   any given linear or nonlinear recovery map. Computing along the above lines
$$
\bar y_w = \argmin_{y\in\cY}\|\cR(A(w),y)\|_{\V'}
$$
we have
\begin{eqnarray}
\label{yest}
&&\|u(y^*)- u(\bar y_w)\|_\U \le  \|u(y^*)- A(w)\|_\U + \|A(w)- u(\bar y_w)\|_\U\nonumber\\
&&\qquad\quad \le E_\wc(A,\cM,\W) + \kappa \dist\,(A(w),\cM)_\U\le (1+\kappa) E_\wc(A,\cM,\W).\qquad\quad
\end{eqnarray}
We consider now the specific elliptic model \eref{Poisson} with affine diffusion coefficients
$a(y)$   given by \eref{affine}. For this model,   it was established in \cite{parest} that for strictly positive $f$ and certain regularity assumptions on $a(y)$ 
as functions of $x\in\Omega$, parameters may be estimated by states. Specifically, when $a(y)\in H^1(\Omega)$ uniformly in $y\in\cY$,
one   has an inverse stability estimate of the form
\be
\label{parest1}
\|a(y)- a(\t y)\|_{L_2(\Omega)}\le C\|u(y)- u(\t y)\|^{1/6}_\U.  
\ee
Thus, whenever the recovery map $A$ satisfies \eref{wishtrue} for some prescribed $\e>0$,
we obtain a parameter estimation bound
of the form
$$
\|a(y^*)- a(\bar y_w)\|_{L_2(\Omega)}\le C \delta_{\kappa\e}(\cM,\W)^{1/6},
$$
Note that when the basis functions $\theta_j$ are $L_2$-orthogonal, $\|a(y^*)- a(\bar y_w)\|_{L_2(\Omega)}$ 
is equivalent to a (weighted) $\ell_2$ norm of $y^*-\bar y_w$.

\subsection{Concluding remarks}

The linear or piecewise linear recovery scheme hinges on the ability to approximate a solution manifold effectively
by linear or affine spaces, globally or locally. As explained earlier this is true for problems of elliptic or parabolic type 
that may include convective terms as long as they are dominated by diffusion. This may however no longer be the case 
when dealing with pure transport equations or models  involving strongly dominating
convection. 

An interesting alternative would then be to adopt a stochastic model  according to \eref{meansquare} and \eref{exp}
that allows one to view the construction of the recovery map as a regression problem. In particular, when dealing with
transport models, a natural candidate for parametrizing a reduced model are {\em deep neural networks}.
However, properly adapting the architecture, regularization and training principles pose wide open questions 
addressed in current work in progress.



\begin{thebibliography}{99}

  
\bibitem{Adcock} B. Adcock, A. C. Hansen, and C. Poon, {\it Beyond Consistent Reconstructions: Optimality and Sharp Bounds for Generalized Sampling, and Application to the Uniform Resampling Problem}, SIAM J. Math. Anal. 45 (2013), 3132-3167.


\bibitem{BC16}
 {M. Bachmayr, A. Cohen},
 {Kolmogorov widths and low-rank approximations of parametric elliptic {PDE}s},
 {Math. Comp.},
 {86}(2017),
 {701--724}. 

\bibitem{BCDDPW1}
 P. Binev, A. Cohen,  W. Dahmen, R. DeVore, G. Petrova, P. Wojtaszczyk, 
 \textit{Convergence Rates for Greedy Algorithms in Reduced
Basis Methods},
SIAM J. Math. Anal., \textbf{43} (2011), 1457--1472.

\bibitem{BCDDPW2}
 P. Binev, A. Cohen,  W. Dahmen, R. DeVore, G. Petrova, P. Wojtaszczyk,  \textit{Data Assimilation in Reduced Modeling}, {SIAM J. Uncertainty Quantification}, \textbf{5}(1) (2017), 1--29.
 
\bibitem{parest}
A. Bonito, A. Cohen, R. DeVore, G. Petrova, and G. Welper, Diffusion Coefficients Estimation for Elliptic Partial Differential Equations,  SIAM J. Math. Anal., 49(2)(2017), 1570--1592. 

\bibitem{BS}
D. Broersen, R. Stevenson, A robust Petrov-Galerkin discretisation of convection-diffusions equations, Comput. Math. Appl. 68(11), 1605--1618 (2014). 
 
  
\bibitem{BDS}
D. Broersen, W. Dahmen, R. Stevenson. On the stability of DPG formulations of transport equations, 
Math. Comp., 87(311)(2018), 1051--1082.



\bibitem{Maday}
A. Buffa, Y. Maday, A.T. Patera, C. Prud'homme, G. Turicini. 
A piori convergence of the greedy algorithm for the parametrized reduced bases,
ESAIM Math. Model. Numer. Anal., 46 (03) (2012), 595-603.

\bibitem{CMN}
  A. Cohen, W. Dahmen, O. Mula and J. Nichols,  {Reduced models for nonlinear state and parameter estimation}, preprint Jan. 2020.

\bibitem{CDD-rand}
A. Cohen, W. Dahmen, R. DeVore,  {Reduced Basis Greedy Selection Using Random Training Sets}, preprint Oct. 2018,
to appear in: ESAIM: Math. Model. Numer. Anal., 
http://arxiv.org/abs/1810.09344 [math.NA].

\bibitem{CDW}
A. Cohen, W. Dahmen, G. Welper, 
Adaptivity and variational stabilization for convection-diffusion equations, 
ESAIM: Math. Model. Numer. Anal., 46(5)(2012), 1247--1273.

\bibitem{CDDFMN}
A. Cohen, W. Dahmen, R. DeVore, J. Fadili, O. Mula,  J. Nichols,
  {Optimal reduced model algorithms for data-based state estimation},
http://arxiv.org/abs/1903.07938

\bibitem{CD} A. Cohen and R. DeVore, {Approximation of high-dimensional PDEs}, 
Acta Numerica,   {24} (2015), 1--159.

\bibitem{CDnw}
A. Cohen, R. DeVore,  Kolmogorov widths under holomorphic mappings, IMA J. Numer. Anal., 36(1)(2016), 1--12.

\bibitem{CDS}
{A. Cohen, R. DeVore, C. Schwab},
    {Convergence rates of best {$N$}-term {G}alerkin approximations
              for a class of elliptic s{PDE}s},
    {Found. Comput. Math.},
      {10}(6) {(2010)}, {615--646},

\bibitem{D}
W. Dahmen, How to Best Sample a Solution Manifold?, 
in: Sampling Theory, a Renaissance,  Applied and Numerical Harmonic Analysis, Series ed. G\"otz E. Pfander, Birkh\"auser, ISBN 978-3-319-19748-7, DOI 10.1007/978-3-319-19749-4.
http://arXiv:1503.00307 [math.NA]
\bibitem{DHSW}
W. Dahmen, C. Huang, C.Schwab, G.Welper, 
Adaptive Petrov-Galerkin methods for first order transport equations, 
 SIAM J. Numer. Anal., 50(5) (2012), 2420--2445.

\bibitem{DPW}
W. Dahmen, C. Plesken, G. Welper, Double greedy algorithms: reduced basis methods for transport dominated problems,
ESAIM: Math. Model. Numer. Anal., 48(3) (2014), 623--663.

\bibitem{DePW}
{R., DeVore,  G. Petrova, P. Wojtaszczyk.}
     {Greedy algorithms for reduced bases in {B}anach spaces},
     {Constr. Approx.},
     {37}(3) {(2013)}, {455--466},

\bibitem{MPPY}
Y. Maday,  A.T. Patera, J.D. Penn and M. Yano, {A parametrized-background data-weak approach to variational data assimilation: Formulation, analysis, and application to acoustics},   {Int. J. Numer.  Meth. Eng.,  Special Issue:Model Reduction}, {102} (5)
(2015), 931--1292

\bibitem{MR}
C. A. Micchelli, T. J. Rivlin,    A survey of optimal recovery, ``Optimal Estimation in Approximation Theory'' (Eds. C. A. Micchelli and T. J. Rivlin), Plenum, N. Y., 1977, 1--54.

\bibitem{Patera}
G. Rozza,  D. B. P.  Huynh,  A.T. Patera. 
      {Reduced basis approximation and a posteriori error estimation
              for affinely parametrized elliptic coercive partial
              differential equations: application to transport and continuum
              mechanics},
    {Arch. Comput. Methods Eng.},
    15(3) (2008),
    {229--275},
       
       \bibitem{SW}
       R. Stevenson, J. Westerdiep, Stability of Galerkin discretizations of a mixed space-time variational formulation of parabolic evolution equations, Preprint, February 2019. Submitted.
       
\bibitem{Stuart}
A. M. Stuart. Inverse problems: A Bayesian perspective. Acta Numerica, 19 (2010), 451--559.

  

\end{thebibliography}
\end{document}